\documentclass[12pt]{article}
\usepackage{diagrams}
\diagramstyle{midshaft,PostScript=dvips,nohug}
\usepackage{amsmath}
\usepackage{amssymb}
\usepackage[T1]{fontenc}
\usepackage{mathptmx}

\newcommand{\T}{\mathbb T}
\newcommand{\A}{\mathbb A}
\newcommand{\E}{\mathcal{E}}

\newcommand{\set}{\mbox{Set}}

\DeclareMathOperator{\Hom }{Hom}
\DeclareMathOperator{\Alg }{-Alg}
\DeclareMathOperator{\incl }{incl}

\newcommand{\pit}{\pitchfork}
\newarrow{Into}  C--->
\newarrow{Implies} ===={=>}
\newtheorem{thm}{Theorem}

\newtheorem{prop}[thm]{Proposition}

\title{Duality for generic algebras}
\author{Anders Kock\\ \small University of Aarhus \normalsize}
\date{}

\begin{document}
\maketitle
\section*{Introduction}For any finitary algebraic theory $\T$, one has a certain topos $\E$ 
with a 
$\T$-algebra object $R$ in it, which classifies $\T$-algebra objects 
in  arbitrary toposes; this $\T$-algebra $R \in \E$ is called the 
{\em generic} $\T$-algebra. The description of this ``classifying 
topos''$\E$, and the $\T$-algebra $R$ in it, is 
simple and well known: $\E$ is the presheaf topos $[FP\T , \set]$, 
where $FP\T$ is the category of finitely presented $\T$-algebras, 
and $R$ is the ``forgetful functor'' $FP\T \to \set$; see e.g.\ 
\cite{MM} Ch.\ VIII, or \cite{JE} Ch.\ D.3.  For any $C\in 
FP\T$, we have two particular objects in $\E$, namely  $y(C)$, where 
$y$ is the Yoneda embedding, and $\gamma^{*}(C)$, where $\gamma^{*}$ 
is left adjoint to the global sections functor $\gamma_{*} : \E \to 
\set$. There is a canonical pairing
$$\gamma^{*}(C) \times y(C) \to R,$$
which we shall describe. The two exponential transposes of this map 
give rise to some duality isomorphisms.
\footnote{ The main result presented here (Theorem \ref{thmax}) 
was presented at the 17th PSSL in 1980 in Sussex, and  
 announced in \cite{GAGD} (1981). I apologize for the long delay in 
publishing a complete account.
}

\section{Generalities on exponential objects}
Let $\E$ ba a Cartesian closed category, and let $Q$ and $R$ be arbitrary objects in it. 
The exponential object $R^{Q}$, which we shall denote $Q\pit R$, is 
characterized by a universal property: it gives rise to the processes of 
``exponential transposition''; these transpositions associate to $k: P\times Q \to R$ 
a map 
$i:Q \to P\pit R$, as well as a map $j: P \to Q \pit R$; the one comes 
about from the other via the symmetry $P\times Q \cong Q\times P$. 
The map $i:Q \to P\pit R$ and $j: P \to Q \pit R$ are {\em 
twisted exponential adjoints} of each other. They are related by a 
map $\delta$, in a commutative diagram
$$\begin{diagram}P&\rTo^{\delta}& (P\pit R)\pit R\\
\dTo^{j}&\ldTo_{i\pit R}&&\\ 
Q\pit R &&&
\end{diagram}$$
where $\delta$ itself is the twisted exponential adjoint of the 
identity map of $P\pit R$. We refer to $\delta$ as a ``Dirac'' map - 
a natural ``embedding'' of an object $P$ into its double dual w.r.to 
a fixed $R$. (It need not be a monic, but often is). All this belongs 
to the elementary theory of Cartesian closed categories, or to ``pure 
lambda calculus''.

 For the case where $\E$ is a presheaf 
topos $[\A ,\set ]$ , we shall recall 
the one of these processes of exponential transposition in elementary 
terms (the other one then comes by the symmetry). The 
$\set$-valued hom-functor of $\A$ we denote by square brackets like 
$[X,Y]$. First, 
we describe the exponential object $Q\pit R$ itself, 
namely for $B\in \A $, 
$$(Q\pit R)(B)= \int _{g\in B/\A} \hom  (Q(X),R(X)),$$
where $X$ denotes the codomain of $g$, and where $\hom$ denotes the 
hom functor for the category of sets. (We shall also use the notation 
$\int_{g:B\to X}$, for $\int _{g\in B/\A}$.) Thus, for an object
$S\in [\A , \set]$ to qualify for the name $Q\pit R$, the object $S$ should for 
each $B \in \A$ be equipped with maps $\pi _{g}: S \to \hom(Q(X), R(X))$ for 
each object $g\in B/\A$, ($X$ denoting the codomain of $g$), subject to certain naturality conditions and 
a certain universal property. Then in terms of the maps $\pi_{g}$, the 
exponential transpose of a map $k:P\times Q \to R$ is given as the 
map $i=\hat{k}: P \to 
Q\pit R$, with $(\hat{k})_{B}$ the unique 
map  such that for  each $g\in B/\A$, we have 
 \begin{equation}\label{transposex}\pi_{g}\circ (\hat{k})_{B}=\widehat{(k_{X})}\circ P(g) 
.\end{equation}
 Here $k_{X}$ is a map $P(X)\times Q(X) \to R(X)$ in the category of 
sets, so its exponential transpose $\widehat{(k_{X})}:P(X) \to \hom 
(Q(X),R(X))$ makes 
immediately sense. 

Two cases will be of particular interest, namely the case where $Q$ is 
re\-pre\-sentable, and where $Q$ is constant. 

In the case where $Q$ is representable, say $Q=y(C)$ for $C\in \A$, 
one has a well known explicit presentation of $Q\pit R$, provided binary 
coproducts exist in $\A$.  Then 
$y(C)\pit R $ may be taken to be $R \circ (-\otimes C): \A \to \set$; 
 let us be explicit about the maps $\pi _{g}$ which qualify
$R \circ (-\otimes C)$ as $y(C)\pit R$. So let $B\in \A$, and let 
$g:B\to X$.
Then 
$$\pi _{g}:R(B\otimes C) \to \hom ([C,X], R(X))= \Pi _{f\in 
[C,X]}R(X)$$ is described by describing its  $f$-coordinate, for $f\in 
[C,X]$: \begin{equation}\label{pixx}p_{f}\circ \pi_{g}:= 
R(\{g,f\})\end{equation}
where $\{g,f\}: B\otimes C \to X$ denotes that map out of the coproduct 
whose components are $g$ and $f$, respectively, and where $p_{f}$ 
denotes the projection to the $f$- factor of the product (or, seeing 
the latter as $\Hom ([C,X],R(X))$, as evaluation at the element $f\in 
[C,X]$).

In the case where $Q$ is constant $Q=\gamma^{*}(C)$ for some set $C$, i.e.\ $Q(X)=C$ 
for all $X\in \A$, we have the following simple presentation of 
$(\gamma^{*}(C)\pit R)(B)$; namely 
\begin{equation}\label{constxx}(\gamma^{*}(C)\pit R)(B)= \Hom 
(C,R(B)), \end{equation} which qualifies for this name 
by virtue of  $\pi_{g} = \Hom (C, R(g))$, for $g:B \to X$. 
This can also be seen from the $\int$ formula for
$(\gamma^{*}(C)\pit R )(B)$;  for $$\int _{g:B\to X}\Hom 
(C,R(X))\cong \Hom (C, 
\int_{g:B\to X}R(X)) \cong \Hom (C, R(B)),$$
using that $\int_{g:B\to X}R(X) \cong R(B)$ via $\pi_{1_{B}}$, by Yoneda's Lemma.
 
\section{$\T$-algebras in a topos}

Let $\E$ be a topos. The 
category $\T \Alg (\E )$ of $\T$-algebras in $\E$ will be monadic 
over $\E$  (\cite{LES} and \cite{JW} Chapter V). This monad will in 
fact be $\E$-enriched; see (the proof of) Lemma 5.5 in \cite{JW}. We denote 
the $\E$ enriched hom functor of it by 
$\pit_{\T}$. If $X$ and $Y$ are objects in $\E$ carrying 
$\T$-structures,  have $X\pit_{T}Y \subseteq X\pit Y.$
If $S\in \E$, and $Y$ carries $\T$-structure, then also $S\pit Y$ 
inherits a $\T$-structure, in a canonical way. But for $T$-algebras 
$X$ and $Y$, the $T$-structure, which $X\pit Y$ inherits from $Y$, need not restrict to one 
on $X\pit_{T}Y$, unless $\T$ is a commutative theory (and e.g.\ the theory of 
commutative rings is not commutative).

If $R\in \T\Alg (\E )$, we have the category $R/(\T\Alg (\E ))$, 
which we call the category of $R$-algebras, denoted $R\Alg (\E)$ or 
just $R\Alg$. It, too, is enriched in $\E$, with enriched hom 
functor denoted $\pit_{R}$. (When $\T$ is the theory of 
commutative rings, and $R$ is such a ring, the terminology 
``$R$-algebra'' agrees with the standard use in commutative algebra; 
however, $\pit_{R}$ may in this context mean something different, 
namely the set (or object)   of ``$R$-linear 
maps'', 
as relevant in the theory of Schwartz distributions, (typically with $R = 
{\mathbb R}$ or $= {\mathbb C}$); see also \cite{PRSFA} and \cite{CMTD}, which deal 
with such cases.)   

When $X$ and $Y$ are $R$-algebras, we have 
$X\pit_{R}Y \subseteq X\pit_{\T}Y \subseteq X\pit Y$. 

We shall use $\otimes$ to denote finite coproducts in the 
category of $\T$-algebras, because one primary example is that of some category of 
commutative rings; and also, because the notation $+$ may incorrectly suggest that 
products $\times$ distribute over the coproduct $+$. There is an 
$\E$-enriched (monadic) adjointness between $\T$-algebras in $\E$, and 
$R$-algebras (in $\E$), whose left adjoint is $R\otimes -$.
\section{The pairing and its transposes}
   
We now specialize to the case where $\A$ is some small category of 
$\T$-algebras, closed under finite coproducts (for instance, $\A$ 
may be the category of finitely presented $\T$-algebras). We consider $\E = [\A 
,\set ]$.   

Let $\gamma_{*}$ be the global-sections  functor of $\E$; it has a 
left adjoint  $\gamma^{*}:\set \to [\A, \set ]= \E $. It associates 
to a set $S$ the functor $\A \to \set$ whose value is the functor 
with constant value $S$. Since $\gamma^{*}$ preserves finite limits, 
it preserves algebraic structure, thus if $C \in \A$, $C$ carries 
$\T$-structure, and therefore 
$\gamma^{*}(C)$ will carry structure of $\T$-algebra in $\E$.

By $R$, we denote the ``forgetful'' functor $\A \to \set$. As an 
object in $\E$, it carries  $\T$-algebra structure, and hence so does 
any object of the form $S\pit R$.

Let $C \in \A$. 
We describe a map 
$k^{C}:\gamma^{*}(C) \times y(C) \to R$, and its two exponential 
transposes $i: y(C) \to 
\gamma^{*}(C)\pit R$ and  $j: \gamma^{*}(C)\to y(C)\pit R$. 
Since $C$ will be fixed in the following, we shall omit the upper 
index $C$ from notation. The map $k=k_{C}$ will be a $\T$-homomorphism in the 
first variable - this notion makes sense for enriched monads, see 
\cite{BCCM}. Therefore, by loc.cit., $i$ will 
factor through $\gamma^{*}(C)\pit_{\T}R \subseteq \gamma^{*}(C) \pit 
R$, and $j$ will be a $\T$-homomorphism.

First, we describe $k$ by describing $k_{B}: \gamma^{*}(C)(B) \times 
y(C)(B) \to R(B)$. Recall that $\gamma^{*}(C)(B)=C$ for all $B$, that
$y(C)(B) = [C,B]$ (where square brackets denote the hom functor of 
$\A$), and recall that $R(B)=B$. Then the map 
$k_{B}: C \times [C,B] \to  B$
is simply the evaluation map $(c,f) \mapsto f(c)$ 
for $c\in C$ and $f:C\to B$ in $\A$. It is a $\T$-homomorphism in the 
variable $c$ because $f$ is a $\T$-homomorphism. Thus, we have three 
maps 
\begin{align}k:& \gamma^{*}(C) \times y(C) \to R\mbox{ , a $\T$-homomorphism in the 
first variable}\\
i:& y(C)\to \gamma^{*}(C)\pit_{\T} R \subseteq \gamma^{*}(C)\pit R\\
j: &\gamma^{*}(C) \to y(C)\pit R\mbox{ , a 
$\T$-homomorphism}\end{align}
Using the explicit 
description of exponential transposition given above, and using 
(\ref{constxx}), it is 
straightforward to see that $i_{B}: y(C)(B) \to (\gamma^{*}(C)\pit 
R)(B) $ is given by the following recipe. We need to give a map 
$i_{B}: 
y(C)(B) \to \Hom  (C, R(B))$
this is just the inclusion 
$[C,B] \subseteq \Hom (C,B)$, which is the $B$-component of the 
inclusion $\gamma^{*}(C)\pit_{\T} R \subseteq \gamma^{*}(C)\pit R $.
  This proves
\begin{prop}\label{yxx}The map $i: y(C) \to \gamma^{*}(C)\pit_{\T}R$ is an 
isomorphism.
\end{prop}

The right hand side of this isomorphism deserves the name 
``$\mbox{Spec}_{R}C$'', since it takes finite colimits to limits, and, 
for $C=$ the free $\T$-algebra in one generator, it   gives $R$, cf.\ 
\cite{SDG} I.12; this notion does not depend on the specifics of the 
topos $\E$.

\medskip

Next, we study the $\T$-homomorphism $j: \gamma^{*}(C) \to y(C)\pit 
R$. For $B\in \A$, we consider
$j_{B}: \gamma^{*}(C)(B) \to (y(C)\pit R)(B)$. Let us first note that 
 the  transpose of the set theoretic map $k_{X}: C \times 
[C,X] \to X$ 
  is the 
``Dirac'' map $\widehat{k_{X}}:C\to \Hom ([C,X],X)$. Next, we utilize  
the 
``coproduct'' description of $y(C)\pit R = R(-\otimes C)$, being an 
exponential object $y(C)\pit R$ by virtue 
of the maps $\pi_{g}$ described in (\ref{pixx}) above. In  terms of 
this, we prove 
\begin{prop}\label{jxx} The map $j: \gamma^{*}(C) \to y(C)\pit R = R( -\otimes 
C)$ has for its $B$-compo\-nent just the inclusion map $\incl _{2}: C \to B\otimes 
C$ into the second compo\-nent of the coproduct. 
\end{prop}
{\bf Proof.} Using the explicit characterization of exponential 
transposition given in (\ref{transposex}), it suffices to see that
$\incl _{2}$ has the property that (for $g:B\to X$), $\pi_{g}\circ 
\incl _{2}= \widehat{k_{X}}$ 
- note that the $P(g)$ occuring in (\ref{transposex}) here is an 
identity map. We analyzed above that $\widehat{k_{X}}$ here is the 
relevant Dirac map
$\delta$, so the task is to prove that the upper triangle in the 
following diagram commutes:
$$
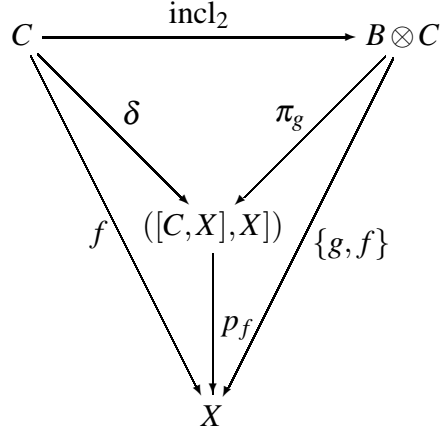
\begin{diagram}C&&\rTo^{\incl _{2}}&& B\otimes C\\
&\rdTo^{\delta}\rdTo(2,4)_{f}&&\ldTo^{\pi_{g}}\ldTo(2,4)_{\{g,f\}}&\\
&&([C,X],X])&&\\
&&\dTo_{p_{f}}&&\\
&&X&&
\end{diagram}$$
and this follows if for all $f:C\to X$, 
it commutes after postcomposition by $p_{f}$, displayed as the 
vertical arrow in the diagram.
Here we write $([C,X],X)$ instead of $\Hom ([C,X],X)$, for typographical 
reasons. The right hand triangle commutes, by construction of 
$\pi_{g}$, and the left hand triangle commutes, by lambda calculus. 
Finally, the outer triangle commutes, by  definition of $\{g,f\}$. 
Therefore, the upper triangle commutes, and this shows that $\incl _{2}$ 
is indeed the claimed exponential transpose. This proves the 
Proposition.

\medskip

We already know from more abstract reasons that $j$ is a 
$\T$-homomorphism; this also appears explicitly from the above 
Proposition, since the coproduct inclusion $C \to B\otimes C$ is a 
$\T$-homomorphism. Now $B\otimes C$ is not only a $\T$-algebra, but 
it is a $B$-algebra by virtue of the coproduct inclusion $\incl _{1}:B \to 
B\otimes C$. Any $\T$-algebra $X$ extends uniquely to a $B$-algebra $B 
\otimes X$, the ``free $B$-algebra in $X$''. 
The canonical extension of the $\T$-algebra morphism
$i_{2}:C\to B\otimes C$ to a $B$-algebra morphism  
$B\otimes C \to B\otimes C$ is clearly the identity map. The 
$R$-algebra structure $R \to y(C)\pit R$ of $y(C)\pit R = (-\otimes C)$ has for its 
$B$-component just  $\incl _{1}$. Since coproducts 
$\otimes$ of $\T$-algebras in a presheaf topos are calculated 
coordinatewise, the free $R$-algebra  $R\otimes \gamma^{*}(C)$ 
on $\gamma^{*}(C)$ has for its $B$ coordinate $B \otimes C$. This 
proves
\begin{thm}\label{jxj}The extension of the $\T$-homomorphism 
$j: \gamma^{*}(C) \to y(C)\pit R$ to an $R$-algebra morphism 
$\overline{j}:R 
\otimes \gamma^{*}(C) \to y(C)\pit R$ is an isomorphism. 
\end{thm}

\noindent {\bf Example 1.} If $\T$ is the theory of commutative rings, 
and $C$ is the ring of dual numbers ${\mathbb Z}[\epsilon]$, then in 
the commutative ring classifier topos $[FP\T , \set ]$, the 
isomorphism $j$ in this Proposition gives in particular the isomorphism of the simplest KL 
axiom, saying that the map $R\times R\to  R^{D}$ is an isomorphism of 
$R$-algebras (with algebra structure on $R\times R$ being ``the ring 
of dual numbers $R[\epsilon ]$). 
(here $D= y({\mathbb Z}[\epsilon]$).

Similarly $R[X]$ (= the free $R$-algebra in one generator) is 
isomorphic, via $j$ for ${\mathbb Z}[X]$, to $R\pit R$.

\medskip

\noindent {\bf Example 2.} If $\T$ is the initial algebraic theory 
(so  $\T$-algebras are just sets), the generic algebra $R$
is called the generic {\em object}, and the classifying topos is called 
the {\em object classifier}, cf.\ \cite{JW} Ch.\ IV (they write $U$ rather 
than $R$). Coproducts $\otimes$ of ``algebras'' are here better denoted $+$; and 
the isomorphism $j$ in this case is a map $ R+1 \to R\pit R$. The ``added'' 
point in the domain of this $j$ is  mapped by $j$ to the identity map of 
$R$.

\medskip

 If $X$ is a $\T$-algebra in $\E$, and $Z$ is an 
$R$-algebra,  we have an isomorphism in $\E$ 
 between $X\pit_{\T}Z$ and $(R\otimes X )\pit_{R}Z$, expressing the 
enrichment of the adjointness between $\T$-algebras in $E$ and 
$R$-algebras (in $\E$).  Using the notion of $R$-algebra and 
this isomorphism, we may reformulate Proposition 
\ref{yxx}. There is no harm in denoting the isomorphism $y(C) \to 
\gamma^{*}C\pit_{\T} R$  of  
Proposition \ref{yxx}  and the isomorphism $y(C) \to 
(R\otimes \gamma^{*}C)\pit_{T} R$, by the same symbol 
$i$: So Proposition \ref{yxx} is reformulated:
\begin{thm}\label{yxxx}
The map $i: y(C) \to (R\otimes \gamma^{*}(C))\pit_{R}R$ is an 
isomorphism.\end{thm}

\section{Duality}

The theme of double dualization occurs in many guises in many areas 
of mathematics. In a Cartesian 
closed category, the simplest is full double dualization functor $(-\pit 
R)\pit R$ into an 
object $R$; there is a natural transformation, whose instantiation at 
$X$ is a map $\delta_{X}:X \to (X\pit R)\pit R$ (where ``$\delta$'' is 
for ``Dirac'', as in Section 1). There are restricted variants of $\delta$, in case $R$ 
carries some algebraic structure, say, of $\T$-algebra. Then one has
$X\to (X\pit R)\pit_{\T}R$ (as studied above); and   in 
case that also $X$ carries $\T$-structure, we have a $\T$-homomorphism $X \to (X\pit_{\T}R)\pit R$, 
obtained by postcomposing $\delta _{X} :X\to (X\pit R)\pit R$ with $s\pit 
R$, where $s$ denotes the inclusion of $X \pit_{\T}R$ into $X \pit 
R$. This composite will also be denoted $\delta_{X}$.
Similarly, if $X$ is an $R$-algebra, we have an $R$-algebra 
homomorphism $\delta_{X}: X \to (X \pit_{R}R)\pit R$.

In our context, the dualization functors (are contravariant and) go 
from ``geometric objects'' (objects in $\E$), to ``algebraic objects'' 
($\T$-algebras, say), and vice versa; the object $R$ is, as a geometric object, the {\em 
line}, but it is canonically endowed with a $\T$-algebra structure, so it lives 
in both worlds. Similarly, $C \in \A$ is a $\T$-algebra, but it 
represents a geometric object $y(C)$. This is the reason for the 
title of the announcement \cite {GAGD}.  

Duality Theorems often have as conclusion that one or the other of 
the Dirac maps mentioned above is an isomorphism. Such  duality results occur in 
our context, as Corollaries of the results above; we shall prove
\begin{thm}\label{thmax} For any $C \in \A$, we have that
$$\delta_{y(C)}: y(C) \to (y(C)\pit R)\pit_{R}R$$
is an isomorphism in $\E$.
\end{thm}
This one may see as a ``Gelfand duality'' result; it  will  
follow from a duality result concerning the $R$-algebra $R\otimes 
\gamma^{*}(C)$:
\begin{thm}\label{thmbx}For any $C \in \A$, we have that
$$\delta_{R\otimes 
\gamma^{*}(C)}: R\otimes \gamma^{*}(C) \to ((R\otimes \gamma^{*}(C) 
)\pit_{R}R)\pit R $$
is an isomorphism of $R$-algebras in $\E$.
\end{thm}

We begin by proving Theorem \ref{thmbx}. We replace the 
pairing $k:\gamma^{*}(C)\times y(C) \to R$ (which is a $\T$-homomorphism in 
the first variable) by its extension to  a 
pairing
\begin{equation}\label{overlx}\overline{k}:
( R\otimes \gamma^{*}(C))\times y(C) \to R,\end{equation}
(which is an $R$-algebra morphism in the first variable), and its two 
exponential transposes $\overline{i}$ and $\overline{j}$; 
here, $\overline{i}$ factors as 
$$\begin{diagram}
y(C)&\rTo^{i}&(R\otimes \gamma^{*}C)\pit_{R}R&\rTo^{ s}&(R\otimes 
\gamma^{*}C)\pit R \end{diagram}$$
where $s$  denotes the inclusion of the $\pit_{R}$ into $\pit$; and 
$\overline{j}$ is the extension of the $\T$-homomorphism $j$ to an 
$R$-homomorphism. 
By ``pure lambda calculus'', as steted in Section 1, we have commutativity of the upper left 
triangle in
$$\begin{diagram}R\otimes \gamma^{*}C&\rTo^{\delta}& ((R\otimes 
\gamma^{*}C)\pit R)\pit R\\
\dTo^{\overline{j}}_{\cong}&\ldTo^{\overline{i}\pit R}&\dTo_{s\pit R}\\ 
yC\pit R &\lTo _{i\pit R}^{\cong}& ((R\otimes 
\gamma^{*}C)\pit_{R}R)\pit R.
\end{diagram}$$ 
The composite $(s\pit R)\circ \delta$ in this diagram is the Dirac 
map considered in the statement of the Theorem. From the commutativity of the 
diagram, and the fact that $\overline{j}$ and $i\pit R$ are 
isomorphisms (Theorems  and \ref{jxj}), we deduce that the $\delta$ of 
the Theorem is an isomorphism, as claimed.

\medskip

To prove Theorem \ref{thmax}, we apply the dualization functor 
$-\pit_{R}R$ to the isomorphism of Theorem \ref{thmbx}, and conclude 
that we get an isomorphism $\delta_{X} 
\pit_{R}R$ (in $\E$) from the 
right to the left in
$$\begin{diagram}(R\otimes \gamma^{*}C)\pit _{R}R & 
\pile{\lTo^{\delta_{X} 
\pit_{R}R}\\ \rTo_{\delta_{Y}}}& 
(((R\otimes \gamma^{*}C)\pit_{R}R)\pit)R)\pit_{R}R
\end{diagram}$$
where $X:=R\otimes \gamma^{*}C$ and $Y:= (R\otimes \gamma^{*}C)\pit 
_{R}R$. However, the map $\delta_{Y}$ here is a splitting of 
$\delta_{X}\pit_{R}R$, by the triangle identity for the adjointness
$$\begin{diagram}\E &\pile{\rTo^{-\pit R}\\ \lTo_{-\pit_{R}R}}& 
(R\mbox{-Alg})^{op}
\end{diagram}$$
(or by ``pure lambda-calculus''). But a splitting of an isomorphism is an 
isomorphism, so we conclude that $\delta_{Y}$ is an isomorphism. Now 
by Theorem \ref{yxxx}, the $Y$ here is isomorphic to $y(C)$, whence 
also $\delta_{y(C)}$ is an isomorphism, proving Theorem \ref{thmax}.

\medskip

The theorems \ref{thmax} and \ref{thmbx} together provide an example 
of a {\em complete pairing} in a sense to be described now. I don't 
many examples, but the notion itself seems to have an aesthetic 
value. Let $T_{1}$ and $T_{2}$ be $\E$-enriched (= strong) monads on a 
Cartesian closed category $\E$, and let $R\in \E$ be an object 
equipped with algebra structures for both the monads; 
 these two structures should commute with each other, in the sense 
described in \cite{DD}, Section 4. Let $P$ be a $T_{1}$-algebra and 
$Q$ a $T_{2}$-algebra, and let $k:P\times Q \to R$ be a map which is 
a $T_{1}$-homomorphism in the first variable, and a 
$T_{2}$-homomorphism in the second varaible. There results, by general 
theory, a $T_{1}$- homomorphism $i:P \to Q\pit_{T_{2}}R$, and a 
$T_{2}$-homomorphism $j: Q 
\to P\pit_{T_{1}}R$. Then $k$ deserves the name complete pairing if 
both $i$ and $j$ are isomorphisms. A complete pairing gives rise to 
two Dirac maps, both of which are isomorphisms, and may in fact be 
described in these terms.

If $T_{1}$ is the monad whose 
algebras are $R$-algebras, as in the theorems quoted, and $T_{2}$ 
is the identity monad, then 
 the $\overline{k}$ considered in (\ref{overlx}) satisfies the conditions.
Another example is with $T_{1}$ the theory of boolean algebras, $T_{2}$ 
the inital theory, $\E$ the category of sets, and $R=2$. Then for any 
finite set $C$, one has a complete pairing, namely the evaluation map 
$(C\pit R)\times C \to R$. This example one may see as the 
origin of Stone duality.

\small \noindent Diagrams were made with Paul Taylor's ``Diagrams'' package.

\bigskip

\noindent
Anders Kock, Dept.\ of Math., University of Aarhus, Denmark

\noindent kock (at) math.au.dk

\noindent December 2014.

\end{document}